\documentclass[a4paper,leqno]{amsart}
\usepackage{amsmath, amssymb, amsthm}

\usepackage[mathscr]{eucal}

\newtheorem{thm}{\textbf{Theorem}}[section]
\newtheorem{prop}[thm]{\textbf{Proposition}}
\newtheorem{lem}[thm]{\textbf{Lemma}}
\newtheorem{cor}[thm]{\textbf{Corollary}}

\theoremstyle{definition}
\newtheorem{defn}[thm]{{\rm Definition}}
\newtheorem{prob}[thm]{{\rm Problem}}

\newcommand{\olapla}{\overline\bigtriangleup}
\newcommand{\onabla}{\overline\nabla}

\newcommand{\fpt}{\frac{\partial}{\partial t}}
\newcommand{\p}{\phi}
\newcommand{\ka}{\kappa}
\newcommand{\g}{\gamma}

\title{The second variational formula of the $k$-energy and $k$-harmonic curves}

\author{Shun Maeta}

\curraddr{Nakakuki 3-10-9 Oyama-shi Tochigi
Japan}
\email{shun.maeta@gmail.com}


\begin{document}




\begin{abstract}
In \cite{jell1}, J.Eells and L. Lemaire introduced $k$-harmonic maps,
 and Wang Shaobo \cite{ws1} showed the first variation formula.
 In this paper, we give the second variation formula of $k$-energy,
 and give a notion of index, nullity and weakly stable. 
 We also study $k$-harmonic maps into the product Riemannian manifold,
 and $k$-harmonic curves into a Riemannian manifold with constant sectional
curvature,
 and show their non-trivial solutions.
\end{abstract}

\maketitle


\vspace{10pt}
\begin{flushleft}
{\large {\bf Introduction}}
\end{flushleft}
Theory of harmonic maps has been applied into various fields in differential geometry.
 The harmonic maps between two Riemannian manifolds are
 critical maps of the energy functional $E(\p)=\frac{1}{2}\int_M\|d\p\|^2v_g$, for smooth maps $\p:M\rightarrow N$.
 
On the other hand, in 1981, J. Eells and L. Lemaire \cite{jell1} proposed the problem to consider the {\em $k$-harmonic maps}:
 they are critical maps of the functional 
 \begin{align*}
 E_{k}(\p)=\int_Me_k(\p)v_g,\ \ (k=1,2,\dotsm),
 \end{align*}
 where $e_k(\p)=\frac{1}{2}\|(d+d^*)^k\p\|^2$ for smooth maps $\p:M\rightarrow N$.
G.Y. Jiang \cite{jg1} studied the first and second variation formulas of the bi-energy $E_2$, 
and critical maps of $E_2$ are called {\em biharmonic maps}. There have been extensive studies on biharmonic maps.
 
In 1989 Wang Shaobo \cite{ws1} studied the first variation formula of the 
$k$-energy $E_k$,
 whose critical maps are called $k$-harmonic maps.
 Harmonic maps are always $k$-harmonic maps by definition.
 In this paper, we study $k$-harmonic maps and show the second variational formula of $E_k$. 
 
In $\S \ref{preliminaries}$, we introduce notation and fundamental formulas of the tension field.

In $\S \ref{k-harmonic}$, we recall $k$-harmonic maps.
 
In $\S \ref{second k}$, we calculate second variation of the $k$-energy $E_k(\p)$.

In $\S\ref{product}$, we show the reduction theorem of $k$-harmonic maps into the product spaces.

Finally, in $\S\ref{constant}$, we study $k$-harmonic curve into Riemannian manifold with constant sectional curvature,
 and get non-trivial solution of $k$-harmonic curve.
 Furthermore, we determine the ODE of the $3$-harmonic curve equation into a sphere.  



\section{Preliminaries}\label{preliminaries}
Let $(M,g)$ be an $m$ dimensional Riemannian manifold,
 $(N,h)$ an $n$ dimensional one,
 and $\p:M\rightarrow N$, a smooth map.
 We use the following notation.
 The second fundamental form $B(\p)$
 of $\p$ is a covariant differentiation $\widetilde\nabla d\p$ of $1$-form $d\p$,
 which is a section of $\odot ^2T^*M\otimes \p^{-1}TN$.
For every $X,Y\in \Gamma (TM)$, let
 \begin{equation}
 \begin{split}
 B(X,Y)
&=(\widetilde\nabla d\p)(X,Y)=(\widetilde\nabla_X d\p)(Y)\\
&=\overline\nabla_Xd\p(Y)-d\p(\nabla_X Y)=\nabla^N_{d\p(X)}d\p(Y)-d\p(\nabla_XY). 
 \end{split}
 \end{equation}
 Here, $\nabla, \nabla^N, \overline \nabla, \widetilde \nabla$ are the induced connections on the bundles $TM$,
 $TN$, $\p^{-1}TN$ and $T^*M\otimes \p^{-1}TN$, respectively.
 
 If $M$ is compact,
 we consider critical maps of the energy functional
 \begin{align}
 E(\p)=\int_M e(\p) v_g,
 \end{align}
where $e(\p)=\frac{1}{2}\|d\p\|^2=\sum^m_{i=1}\frac{1}{2}\langle d\p(e_i),d\p(e_i)\rangle$
 which is called the {\em enegy density} of $\p$, and the inner product 
 $\langle \cdot ,\cdot \rangle$ is a Riemannian metric $h$. 
 The {\em tension \ field} $\tau(\p)$ of $\p$ is defined by
 \begin{align}
 \tau(\p)=\sum^{m}_{i=1}(\widetilde \nabla d\p)(e_i,e_i)=\sum^m_{i=1}(\widetilde \nabla _{e_i}d\p)(e_i).
 \end{align}
 Then, $\p$ is a {\em harmonic map} if $\tau(\p)=0$.
 
 The curvature tensor field $ R^N(\cdot, \cdot)$ of the Riemannian metric on the bundle 
$TN$ is defined as follows :
 \begin{align}
R^N(X,Y)
=\nabla^N_X \nabla^N_Y - \nabla^N _Y \nabla^N_X-\nabla^N_{[X,Y]},
\ \ \ \ (X,Y\in \Gamma (TN)).
\end{align}
 

$\olapla
=\onabla^* \onabla
=-\sum^m_{k=1}(\onabla_{e_k}\onabla_{e_k}
-\onabla_{\nabla_{e_k}e_k}),$ 
is the {\em rough Laplacian}.



\vspace{10pt}

\section{$k$-harmonic maps}\label{k-harmonic}

J. Eells and L. Lemaire \cite{jell1} proposed the notation of $k$-harmonic maps. 
The Euler-Lagrange equations for the $k$-harmonic maps was shown by
Wang Shaobo \cite{ws1}.
In this section, we recall $k$-harmonic maps.

We consider a smooth variation $\{\p_{t}\}_{t\in I_{\epsilon}} (I_{\epsilon}=(-\epsilon, \epsilon))$ of $\p$ with parameters $t$,
 i.e. we consider the smooth map $F$  given by
$$F : I_{\epsilon}\times M\rightarrow N, F (t,p)=\p_{t}(p),$$
where $F (0,p)=\p_{0}(p)=\p(p), $ forall $p\in M$.

The corresponding variational vector field $V$ is given by
\begin{align*}
V(p)=\frac{d}{dt}|_{t=0}\p_{t,0}\in T_{\p(p)}N,
\end{align*}
$V$ are section of $\p^{-1}TN$, i.e. $V\in \Gamma (\p^{-1}TN)$.

\begin{defn}[\cite{jell1}]
For $k=1,2,\dotsm$ the {\em $k$-energy functional}
 is defined by 
\begin{align*}
E_k(\phi)=\frac{1}{2}\int_M\|(d+d^* )^k\phi\|^2v_g,\ \ \phi\in C^{\infty}(M,N).
\end{align*}
Then, $\phi$ is {\em $k$-harmonic} if it is a critical point of $E_k,$ i.e., for all smooth variation $\{\phi_t\}$ of $\phi$ with $\phi_0=\phi$,
\begin{align*}
\left.\frac{d}{dt}\right|_{t=0}E_k(\phi_t)=0.
\end{align*}
We say for a $k$-harmonic map to be {\em proper} if it is not harmonic. 
\end{defn}


\vspace{10pt}

\begin{lem}\label{laps-1}
\begin{align*}
\onabla_{\frac{\partial}{\partial t}} \olapla^{s-1} \tau(F) |_{t=0}
=&-\olapla^s V
+\olapla^{s-1}R^N(V,d\p(e_j) )d\p(e_j)\\
&+\sum^{s-1}_{l=1}\olapla^{l-1}
\{
-\onabla_{e_j}R^N(V,d\p(e_j) )\olapla^{s-l-1}\tau(\p) \\
&\hspace{52pt}-R^N(V,d\p(e_j) )\onabla_{e_j}\olapla^{s-l-1}\tau(\p)\\
&\hspace{52pt}+R^N(V,d\p(\nabla_{e_j}e_j )\olapla^{s-l-1}\tau(\p)
\}.
\end{align*}
\end{lem}

\begin{proof}
For all $\omega \in \Gamma(\phi^{-1}TN)$,
\begin{align*}
\onabla_{\frac{\partial}{\partial t}} \olapla \omega 
&=-\{
\onabla_{\frac{\partial}{\partial t}}
(
\onabla_{e_j}\onabla_{e_j}-\onabla_{\nabla_{e_j}e_j}
)\omega
\}\\
&=-\{
\onabla_{e_j}\onabla_{\frac{\partial}{\partial t}}(\onabla_{e_j}\omega)
+R^N(dF(\frac{\partial}{\partial t}),dF(e_j))\onabla_{e_j}\omega \\
&\hspace{10pt}-\onabla _{\nabla_{e_j}e_j}\onabla_{\fpt}\omega
-R^N(dF(\fpt),dF(\nabla_{e_j}e_j))\omega
\}\\
&=-\{\onabla_{e_j}( \onabla_{e_j}\onabla_{\fpt}\omega+R^N(dF(\fpt),dF(e_j))\omega )\\
&\ \ \ +R^N(dF(\frac{\partial}{\partial t}),dF(e_j))\onabla_{e_j}\omega \\
&\ \ \ -\onabla _{\nabla_{e_j}e_j}\onabla_{\fpt}\omega
-R^N(dF(\fpt),dF(\nabla_{e_j}e_j))\omega
\}.
\end{align*}
Repeating this and using 
\begin{align*}
\onabla_{\fpt}\tau(F)|_{t=0}
=-\olapla V+R^N(V,d\p(e_j) )d\p(e_j),
\end{align*}
we have the lemma. 

\end{proof}

\begin{lem}\label{nablaps-1}
\begin{align*}
\onabla_{\frac{\partial}{\partial t}} \onabla_{e_i} \olapla^{s-1} \tau(F) |_{t=0}
=&-\onabla_{e_i}\olapla^s V
+\onabla_{e_i} \olapla^{s-1}R^N(V,d\p(e_j) )d\p(e_j)\\
&+\sum^{s-1}_{l=1}\onabla_{e_i} \olapla^{l-1}
\{
-\onabla_{e_j}R^N(V,d\p(e_j) )\olapla^{s-l-1}\tau(\p) \\
&\hspace{68pt}-R^N(V,d\p(e_j) )\onabla_{e_j}\olapla^{s-l-1}\tau(\p)\\
&\hspace{68pt}+R^N(V,d\p(\nabla_{e_j}e_j )\olapla^{s-l-1}\tau(\p)
\}\\
&+R^N(V, d\p(e_i))\olapla^{s-1}\tau(\p).
\end{align*}
\end{lem}

\begin{proof}
\begin{align*}
\onabla_{\frac{\partial}{\partial t}} \onabla_{e_i} \olapla^{s-1} \tau(F)
=&\onabla_{e_i}\onabla_{\fpt} \olapla^{s-1}\tau(F)
+R^N(dF(\fpt),dF(e_i))\olapla^{s-1}\tau(F).
\end{align*}
And using Lemma $\ref{laps-1}$, we have the lemma.
\end{proof}

\begin{lem}\label{R}
\begin{align*}
\int_M\langle \onabla_{e_j}R^N(V, d\p(e_j) )V_1-R^N(V, d\p (\nabla_{e_j}e_j)),V_1, V_2 \rangle v_g\\
=-\int_M\langle R^N(V, d\p(e_j) )V_1, \onabla_{e_j}V_2 \rangle v_g,
\end{align*}
$V_1,V_2 \in \Gamma (\p^{-1}TN).$
\end{lem}

\begin{proof}
\begin{align*}
{\rm div} (\langle R^N(V,d\p(e_i))V_1, V_2 \rangle e_i)
=&\langle \nabla_{e_j} \langle R^N(V,d\p(e_i))V_1, V_2 \rangle e_i , e_j\rangle\\
=&\langle \langle \onabla_{e_j} R^N(V, d\p (e_i))V_1,V_2\rangle e_i \\
&+\langle R^N(V,d\p(e_i) )V_1, \onabla_{e_j}V_2\rangle e_i\\
&+\langle R^N (V, d\p (e_i))V_1, V_2 \rangle \nabla_{e_j}e_i, e_j \rangle.
\end{align*}
By Green's theorem, we have
\begin{align*}
0=&\int_M \rm{div} \langle R^N(V,d\p(e_i))V_1, V_2 \rangle e_i v_g\\
=&\int_M \langle \onabla _{e_j}R^N(V,d\p (e_i)V_1,V_2)\rangle \delta_{ij}\\
&+\langle R^N(V, d\p (e_i )V_1, \onabla_{e_j}V_2)\delta_{ij}\\
&+\langle R^N (V, d\p(e_i) )V_1, V_2\rangle \langle \nabla_{e_j}e_i,e_j\rangle v_g.
\end{align*}
Here,
\begin{align*}
\langle R^N (V, d\p(e_i) )V_1, V_2\rangle \langle \nabla_{e_j}e_i,e_j\rangle
=&\langle R^N (V, d\p(\langle \nabla_{e_j}e_i,e_j\rangle e_i) )V_1, V_2\rangle\\
=&-\langle R^N (V, d\p(\nabla_{e_i}e_i) )V_1, V_2\rangle.
\end{align*}
Therefore, we have the lamma.
\end{proof}

\begin{thm}\label{2s-harmonic}
Let $k=2s\ \ \ (s=1,2,\cdots)$,
\begin{align*}
\frac{1}{2}\frac{d}{dt}E_{2s}(\p_t)|_{t=0}=\int_M\langle \tau_{2s}(\p), V\rangle,
\end{align*}
where, 
\begin{align*}
\tau_{2s}(\p)
=&-\olapla^{2s-1}\tau(\p)+R^N(\olapla^{2s-2}\tau(\p),d\p(e_j))e\p(e_j)\\
&+\sum^{s-1}_{l=1}
\{
R^N(\onabla_{e_j}\olapla^{s+l-2}\tau(\p),\olapla^{s-l-1}\tau(\p) )d\p(e_j) \\
&\hspace{32pt}-R^N(\olapla^{s+l-2}\tau(\p),\onabla_{e_i}\olapla^{s-l-1}\tau(\p) )d\p(e_j)
\},
\end{align*}
where, $\olapla^{-1}=0$.
\end{thm}

\begin{proof}
\begin{align*}
E_{2s}(\p)
=&\int_M\langle \underbrace{(d^*d)\cdots (d^*d)}_s\p, \underbrace{(d^*d)\cdots (d^*d)}_s\p\rangle v_g
=&\int_M\langle \olapla^{s-1}\tau(\p), \olapla^{s-1}\tau(\p)\rangle v_g
\end{align*}
By using Lemma $\ref{laps-1}$ and Lemma $\ref{R}$, we calculate $\frac{1}{2}\frac{d}{dt} E_{2s}(\p_t)$,
\begin{equation}\label{s}
\begin{split}
\frac{1}{2}&\frac{d}{dt} E_{2s}(\p_t)|_{t=0}\\
=&\int_M\langle \onabla_{\fpt}\olapla^{s-1}\tau(F), \olapla^{s-1}\tau(F)\rangle v_g |_{t=0}\\
=&\int_M\langle -\olapla^s V+\olapla^{s-1}R^N(V,d\p(e_j) )d\p(e_j)\\
&+\sum^{s-1}_{l=1}\olapla^{l-1}
\{
-\onabla_{e_j}R^N(V,d\p(e_j) )\olapla^{s-l-1}\tau(\p) \\
&\hspace{52pt}-R^N(V,d\p(e_j) )\onabla_{e_j}\olapla^{s-l-1}\tau(\p)\\
&\hspace{52pt}+R^N(V,d\p(\nabla_{e_j}e_j )\olapla^{s-l-1}\tau(\p)
,\olapla^{s-1}\tau(\p)\}
\rangle\\
=&\int_M \langle V, -\olapla^{2s-1}\tau(\p)\rangle v_g\\
&+\int_M \langle V, R^N(\olapla^{2s-2}\tau(\p),d\p(e_j ) )d\p(e_j )\rangle v_g\\
&+\sum^{s-1}_{l=1}\int_M\langle 
-\onabla_{e_j }R^N(V,d\p(e_j ) )\olapla^{s-l-1}\tau(\p) \\
&\hspace{52pt}-R^N(V,d\p(e_j ) )\onabla_{e_j}\olapla^{s-l-1}\tau(\p)\\
&\hspace{52pt}+R^N(V,d\p(\nabla_{e_j}e_j )\olapla^{s-l-1} \tau(\p)
, \olapla^{s+l-2} \tau(\p)\rangle v_g\\
=&\int_M \langle V, -\olapla^{2s-1}\tau(\p)\rangle v_g\\
&+\int_M \langle V, R^N(\olapla^{2s-2}\tau(\p),d\p(e_j ) )d\p(e_j )\rangle v_g\\
&+\sum^{s-1}_{l=1} \{\int_M\langle 
R^N(V,d\p(e_j ) )\olapla^{s-l-1}\tau(\p),\onabla_{e_j}\olapla^{s+l-2}\tau(\p)\rangle v_g \\
&+\int_M \langle -R^N(V,d\p(e_j ) )\onabla_{e_j}\olapla^{s-l-1}\tau(\p), \olapla^{s+l-2}\tau(\p)\rangle v_g\}\\
=&\int_M \langle V, -\olapla^{2s-1}\tau(\p)\rangle v_g\\
&+\int_M \langle V, R^N(\olapla^{2s-2}\tau(\p),d\p(e_j ) )d\p(e_j )\rangle v_g\\
&+\sum^{s-1}_{l=1} \{\int_M\langle 
R^N(\onabla_{e_j}\olapla^{s+l-2}\tau(\p) ,\olapla^{s-l-1}\tau(\p)  )d\p(e_j), V\rangle v_g \\
&-\int_M \langle R^N(\olapla^{s+l-2}\tau(\p) ,\onabla_{e_j}\olapla^{s-l-1}\tau(\p) )d\p(e_j ) , V \rangle v_g\}\\
=&\int_M \langle V, -\olapla^{2s-1}\tau(\p)
+R^N(\olapla^{2s-2}\tau(\p),d\p(e_j ) )d\p(e_j )\\
&+\sum^{s-1}_{l=1} \{
R^N(\onabla_{e_j}\olapla^{s+l-2}\tau(\p) ,\olapla^{s-l-1}\tau(\p)  )d\p(e_j)\\
&-R^N(\olapla^{s+l-2}\tau(\p) ,\onabla_{e_j}\olapla^{s-l-1}\tau(\p) )d\p(e_j )\} \rangle v_g.
\end{split}
\end{equation}

So we have the theorem.
\end{proof}

\begin{thm}\label{2s+1-harmonic}
Let $k=2s+1\ \ \ (s=0,1,2,\cdots),$
\begin{align*}
\frac{1}{2}\frac{d}{dt}E_{2s+1}(\p_t)|_{t=0}=\int_M\langle \tau_{2s+1}(\p), V\rangle,
\end{align*}
where, 
\begin{align*}
\tau_{2s+1}(\p)
=&-\olapla^{2s}\tau(\p)+R^N(\olapla^{2s-1}\tau(\p),d\p(e_j))d\p(e_j)\\
&+\sum^{s-1}_{l=1}
\{
R^N(\onabla_{e_j}\olapla^{s+l-1}\tau(\p),\olapla^{s-l-1}\tau(\p) )d\p(e_j) \\
&\hspace{32pt}-R^N(\olapla^{s+l-1}\tau(\p),\onabla_{e_j}\olapla^{s-l-1}\tau(\p) )d\p(e_j)
\}\\
&+R^N(\onabla_{e_i}\olapla^{s-1}\tau(\p),\olapla^{s-1}\tau(\p))d\p(e_i),
\end{align*}
where, $\olapla^{-1}=0.$
\end{thm}

\begin{proof}
When $s=0$, it is well known harmonic map. So we consider the case of $s=1,2,\cdots $.  
\begin{align*}
E_{2s+1}(\p)
=&\int_M\langle d\underbrace{(d^*d)\cdots (d^*d)}_{s}\p, d\underbrace{(d^*d)\cdots (d^*d)}_{s}\p\rangle v_g\\
=&\int_M\langle \onabla_{e_i}\olapla^{s-1}\tau(\p_t), \onabla_{e_i}\olapla^{s-1}\tau(\p_t)\rangle v_g.
\end{align*}
By using Lemma $\ref{nablaps-1}$ and Lemma $\ref{R}$, we calculate $\frac{1}{2}\frac{d}{dt} E_{2s+1}(\p_t)$,
\begin{align*}
\frac{1}{2}&\frac{d}{dt} E_{2s+1}(\p_t)|_{t=0}\\
=&\int_M\langle \onabla_{\fpt}\onabla_{e_i}\olapla^{s-1}\tau(F), \onabla_{e_i}\olapla^{s-1}\tau(F)\rangle v_g |_{t=0}\\
=&\int_M\langle -\onabla_{e_i}\olapla^s V+\onabla_{e_i}\olapla^{s-1}R^N(V,d\p(e_j) )d\p(e_j)\\
&+\sum^{s-1}_{l=1}\onabla_{e_i}\olapla^{l-1}
\{
-\onabla_{e_j}R^N(V,d\p(e_j) )\olapla^{s-l-1}\tau(\p) \\
&\hspace{52pt}-R^N(V,d\p(e_j) )\onabla_{e_i}\olapla^{s-l-1}\tau(\p)\\
&\hspace{52pt}+R^N(V,d\p(\nabla_{e_j}e_j )\olapla^{s-l-1}\tau(\p)\\
&\hspace{52pt}+R^N(V,d\p(e_i))\olapla^{s-1}\tau(\p)
,\onabla_{e_i}\olapla^{s-1}\tau(\p)\}
\rangle v_g.
\end{align*}
Here, using 
\begin{align*}
\int_M\langle \onabla_{e_i}\omega_1, \onabla_{e_i}\omega_2\rangle v_g=\int_M \langle \olapla \omega_1, \omega_2 \rangle v_g,\\
\end{align*}
where, $\omega_1, \omega_2 \in \Gamma (\p^{-1}TN)$, we have
\begin{equation}\label{s+1}
\begin{split}
\frac{1}{2}&\frac{d}{dt} E_{2s+1}(\p_t)|_{t=0}\\
=&\int_M \langle V, -\olapla^{2s}\tau(\p)\rangle v_g\\
&+\int_M \langle  R^N( V,d\p(e_j ) )d\p(e_j ) ,\olapla^{2s-1}\tau(\p) \rangle v_g\\
&+\sum^{s-1}_{l=1}\int_M\langle 
-\onabla_{e_j }R^N(V,d\p(e_j ) )\olapla^{s-l-1}\tau(\p) \\
&\hspace{52pt}-R^N(V,d\p(e_j ) )\onabla_{e_j}\olapla^{s-l-1}\tau(\p)\\
&\hspace{52pt}+R^N(V,d\p(\nabla_{e_j}e_j )\olapla^{s-l-1} \tau(\p)
, \olapla^{s+l-1} \tau(\p)\rangle v_g\\
&+\int_M\langle R^N(V, d\p(e_i ) )\olapla^{s-1}\tau(\p), \onabla_{e_i}\olapla^{s-1}\tau(\p))\rangle v_g\\
=&\int_M \langle V, -\olapla^{2s}\tau(\p)
+R^N(\olapla^{2s-1}\tau(\p),d\p(e_j ) )d\p(e_j )\\
&+\sum^{s-1}_{l=1} \{
R^N(\onabla_{e_j}\olapla^{s+l-1}\tau(\p) ,\olapla^{s-l-1}\tau(\p)  )d\p(e_j)\\
&\hspace{30pt}-R^N(\olapla^{s+l-1}\tau(\p) ,\onabla_{e_j}\olapla^{s-l-1}\tau(\p) )d\p(e_j )\} \\
&+R^N(\onabla_{e_i}\olapla^{s-1}\tau(\p) , \olapla^{s-1}\tau(\p)  )d\p(e_i ) \rangle v_g.
\end{split}
\end{equation}
So we have the theorem.
\end{proof}

By Theorem $\ref{2s-harmonic}$, $\ref{2s+1-harmonic}$, we have the following \cite{ws1}.
\begin{cor}\label{harmonick-harmonic}
harmonic map is always $k$-harmonic map $(k=1,2,\cdots)$.
\end{cor}

\vspace{10pt}

For $\olapla^{l}$ $(k=1,2,\cdots )$, we have Theorem $\ref{olapla harmonic}$.
 First, we show the following two lemmas.


\begin{lem}\label{olapla1}
Let $l=1,2,\dotsm$. If \ $\onabla_{e_i}\olapla^{(l-1)} \tau(\p)=0\ (i=1,\dotsm,m)$, then $$\olapla^l \tau(\p)=0.$$
\end{lem}

\vspace{10pt}

\begin{proof}
Indeed, we can define a global vector field $X_{\p}\in \Gamma(TM)$ defined by
\begin{align}
X_{\phi}=\sum^m_{j=1}\langle-\overline\nabla_{e_j}\overline\bigtriangleup^{(l-1)}\tau(\phi),
\overline\bigtriangleup^l\tau(\phi)\rangle e_j.
\end{align}
 Then, the divergence of $X_{\p}$ is given as 
\begin{align*}
{\rm div(X_{\phi})}
&=\langle\overline\bigtriangleup^l\tau(\phi),\overline\bigtriangleup^l\tau(\phi)\rangle
+\sum^m_{j=1}\langle-\overline\nabla_{e_j}\overline\bigtriangleup^{(l-1)}\tau(\phi),
\overline\nabla_{e_j}\overline\bigtriangleup^{l}\tau(\phi)\rangle\\
&=\langle\overline\bigtriangleup^l\tau(\phi),\overline\bigtriangleup^l\tau(\phi)\rangle,
\end{align*}
by the assumption. Integrating this over $M$, we have
$$0=\int_M{\rm div}(X_{\phi})v_g
=\int_M\langle\overline\bigtriangleup^l\tau(\phi),\overline\bigtriangleup^l\tau(\phi)\rangle v_g,$$
which implies $\overline\bigtriangleup^l\tau(\phi)=0.$
\end{proof}

\vspace{10pt}

\begin{lem}\label{olapla2}
Let $l=1,2,\dotsm$. If\  $\olapla^l \tau(\p)=0$, then $$\onabla_{e_i}\olapla^{(l-1)}\tau(\p)=0,\ \ (i=1,\dotsm, m).$$
\end{lem}


\begin{proof}
Indeed, by computing the Laplacian of the $2l$-energy density $e_{2l}(\phi)$, we have
\begin{equation}\label{non.2}
\begin{split}
\bigtriangleup e_{2l}(\phi)
=&\sum^m_{i=1}\left\langle \overline \nabla_{e_i}\overline\bigtriangleup ^{(l-1)}\tau(\phi),
\overline \nabla_{e_i}\overline\bigtriangleup ^{(l-1)}\tau(\phi)\right\rangle\\
&\hspace{30pt}
-\left\langle \overline\nabla^*\overline \nabla(\overline\bigtriangleup ^{(l-1)}\tau(\phi)),\overline\bigtriangleup ^{(l-1)}\tau(\phi)\right\rangle \\
=&\sum^m_{i=1}\left\langle \overline \nabla_{e_i}\overline\bigtriangleup ^{(l-1)}\tau(\phi),
\overline \nabla_{e_i}\overline\bigtriangleup ^{(l-1)}\tau(\phi)\right\rangle 
\geq 0.
\end{split}
\end{equation}
By Green's theorem $\int_M\bigtriangleup e_{2l}(\phi) v_g=0,$ and $(\ref{non.2}),$ we have 
$\bigtriangleup e_{2l}(\phi)=0.$
 Again, by $(\ref{non.2}),$ we have
$$\overline \nabla_{e_i}\overline\bigtriangleup ^{(l-1)}\tau(\phi)=0,\ \ \ (i=1,\dotsm,m,\ \ \ l=1,2,\dotsm).$$
\end{proof}

\vspace{3pt}

\begin{thm}\label{olapla harmonic}
Let $l=1,2,\dotsm $. If $\olapla^l\tau(\p)=0$ or $\onabla_{e_i}\olapla^{(l-1)}\tau(\p)=0,\ 
\newline(i=1,2,\dotsm,m)$, then
$\phi:M\rightarrow N$ from a compact Riemannian manifold into a Riemannian manifold is a harmonic map.
\end{thm}


\begin{proof}
By using Lemma \ref{olapla1}, \ref{olapla2}, we have Theorem \ref{olapla harmonic}.
\end{proof}

\vspace{20pt}

\section{The second variational formula of the $k$-energy}\label{second k}

In this section we calculate the second variation of the $k$-energy.

Now let $\p :(M,g)\rightarrow (N,h)$ be a $k$-harmonic map $(k=1,2,\cdots)$.
We consider a smooth variation $\{\p_{t,r}\}_{t,r\in I_{\epsilon}} (I_{\epsilon}=(-\epsilon, \epsilon))$ of $\p$ with two parameters $t$ and $r$, i.e.
 we consider the smooth map $F$  given by
$$F : I_{\epsilon}\times I_{\epsilon}\times M\rightarrow N, F (t,r,p)=\p_{t,r}(p),$$
where $F (0,0,p)=\p_{0,0}(p)=\p(p), $ for all $p\in M$.

The corresponding variational vector field $V$ and $W$ are given by
\begin{align*}
V(p)=\frac{d}{dt}|_{t=0}\p_{t,0}\in T_{\p(p)}N,\\
W(p)=\frac{d}{dr}|_{r=0}\p_{0,r}\in T_{\p(p)}N.
\end{align*}
$V$ and $W$ are section of $\p^{-1}TN$, i.e. $V,W\in \Gamma (\p^{-1}TN)$.

The {\em Hessian} of $E_{k}$ at its critical point $\p$ is defined by
\begin{align*}
H(E_k)_{\p}(V,W)=\frac{\partial ^2}{\partial t\partial r}|_{(t,r)=(0,0)}E_k(\p_{t,r}).
\end{align*}

\begin{thm}
Let $\p: (M,g)\rightarrow (N,h)$ be a $2s$-harmonic map $(s=1,2,\cdots)$. Then, the 
 {\em Hessian} of the $2s-$energy $E_{2s}$ at $\p$ is given by 
\begin{align*}
H(E_{2s})_{\p}(V,W)=\int_M\langle V, J_{2s}(W)\rangle v_g,
\end{align*}  
where,
\begin{align*}
J_{2s}(W)=-I_{2s}+II_{2s}+III_{2s}+IV_{2s}.
\end{align*}
where,
\begin{align*}
I_{2s}=&-\olapla^{2s} W
+\olapla^{2s-1}R^N(W,d\p(e_j) )d\p(e_j)\\
&+\sum^{2s-1}_{l=1}\olapla^{l-1}
\{
-\onabla_{e_j}R^N(W,d\p(e_j) )\olapla^{2s-l-1}\tau(\p) \\
&\hspace{52pt}-R^N(W,d\p(e_j) )\onabla_{e_j}\olapla^{2s-l-1}\tau(\p)\\
&\hspace{52pt}+R^N(W,d\p(\nabla_{e_j}e_j )\olapla^{2s-l-1}\tau(\p)
\},
\end{align*}

\begin{align*}
II_{2s}=&-(\nabla^N_{\olapla^{2s-2}\tau(\p)}R^N)(d\p(e_i),W)d\p(e_i)\\
&-(\nabla^N_{d\p(e_i)}R^N)(W,\olapla^{2s-2}\tau(\p))d\p(e_i)\\
&+R^N(-\olapla^{2s-1} W
+\olapla^{2s-2}R^N(W,d\p(e_j) )d\p(e_j)\\
&\hspace{30pt}+\sum^{2s-2}_{l_2=1}\{\olapla^{l_2-1}
\{
-\onabla_{e_j}R^N(W,d\p(e_j) )\olapla^{2s-l_2-2}\tau(\p) \\
&\hspace{52pt}-R^N(W,d\p(e_j) )\onabla_{e_j}\olapla^{2s-l_2-2}\tau(\p)\\
&\hspace{52pt}+R^N(W,d\p(\nabla_{e_j}e_j )\olapla^{2s-l_2-2}\tau(\p)
\}\}, d\p(e_i) )d\p(e_i)\\
&+R^N(\olapla^{2s-2}\tau(\p), \onabla_{e_i}W)d\p(e_i)\\
&+R^N(\olapla^{2s-2}\tau(\p), d\p(e_i) ) \onabla_{e_i}W,
\end{align*}

\begin{align*}
III_{2s}=&
-(\nabla^N_{\onabla_{e_i}\olapla^{s+l-2}\tau(\p)}R^N)(\olapla^{s-l-1}\tau(\p),W)d\p(e_i)\\
&-(\nabla^N_{\olapla^{s-l-1}\tau(\p)}R^N)(W,\onabla_{e_i}\olapla^{s+l-2}\tau(\p))d\p(e_i)\\
&+R^N(-\onabla_{e_i}\olapla^{s+l-1} W
+\onabla_{e_i}\olapla^{s+l-2}R^N(W,d\p(e_j) )d\p(e_j)\\
&\hspace{30pt}+\sum^{s+l-2}_{l_2=1}\{\onabla_{e_i}\olapla^{l_2-1}
\{
-\onabla_{e_j}R^N(W,d\p(e_j) )\olapla^{s+l-2-l_2}\tau(\p) \\
&\hspace{115pt}-R^N(W,d\p(e_j) )\onabla_{e_j}\olapla^{s+l-2-l_2}\tau(\p)\\
&\hspace{115pt}+R^N(W,d\p(\nabla_{e_j}e_j )\olapla^{s+l-2-l_2}\tau(\p)
\}\}\\
&+R^N(W,d\p(e_i))\olapla^{s+l-2}\tau(\p),\olapla^{s-l-1}\tau(\p))d\p(e_i)\\
&+R^N(\onabla_{e_i}\olapla^{s+l-2}\tau(\p), 
\{-\olapla^{s-l} W
+\olapla^{s-l-1}R^N(W,d\p(e_j) )d\p(e_j)\\
&\hspace{30pt}+\sum^{s-l-1}_{l_2=1}\{\olapla^{l_2-1}
\{
-\onabla_{e_j}R^N(W,d\p(e_j) )\olapla^{s-l-1-l_2}\tau(\p) \\
&\hspace{52pt}-R^N(W,d\p(e_j) )\onabla_{e_j}\olapla^{s-l-1-l_2}\tau(\p)\\
&\hspace{52pt}+R^N(W,d\p(\nabla_{e_j}e_j )\olapla^{s-l-1-l_2}\tau(\p)
\}
\})d\p(e_i)\\
&+R^N(\onabla_{e_i}\olapla^{s+l-2}\tau(\p), \olapla^{s-l-1}\tau(\p) ) \onabla_{e_i}W,\\
\end{align*}
\begin{align*}
IV_{2s}=&
-(\nabla^N_{\olapla^{s+l-2}\tau(\p)}R^N)(\onabla_{e_i}\olapla^{s-l-1}\tau(\p),W)d\p(e_i)\\
&-(\nabla^N_{\onabla_{e_i}\olapla^{s-l-1}\tau(\p)}R^N)(W,\olapla^{s+l-2}\tau(\p))d\p(e_i)\\
&+R^N(-\olapla^{s+l-1} W
+\olapla^{s+l-2}R^N(W,d\p(e_j) )d\p(e_j)\\
&\hspace{30pt}+\sum^{s+l-2}_{l_2=1}\{\olapla^{l_2-1}
\{
-\onabla_{e_j}R^N(W,d\p(e_j) )\olapla^{s+l-2-l_2}\tau(\p) \\
&\hspace{52pt}-R^N(W,d\p(e_j) )\onabla_{e_j}\olapla^{s+l-2-l_2}\tau(\p)\\
&\hspace{52pt}+R^N(W,d\p(\nabla_{e_j}e_j )\olapla^{s+l-2-l_2}\tau(\p)
\},\onabla_{e_i}\olapla^{s-l-1}\tau(\p) )d\p(e_i)\\
&+R^N(\olapla^{s+l-2}\tau(\p), 
\{-\onabla_{e_i}\olapla^{s-l} W
+\onabla_{e_i}\olapla^{s-l-1}R^N(W,d\p(e_j) )d\p(e_j)\\
&\hspace{30pt}+\sum^{s-l-1}_{l_2=1}\{\onabla_{e_i}\olapla^{l_2-1}
\{
-\onabla_{e_j}R^N(W,d\p(e_j) )\olapla^{s-l-1-l_2}\tau(\p) \\
&\hspace{115pt}-R^N(W,d\p(e_j) )\onabla_{e_j}\olapla^{s-l-1-l_2}\tau(\p)\\
&\hspace{115pt}+R^N(W,d\p(\nabla_{e_j}e_j )\olapla^{s-l-1-l_2}\tau(\p)
\}\\
&\hspace{30pt}+R^N(W,d\p(e_i))\olapla^{s-l-1}\tau(\p)\})d\p(e_i)\\
&+R^N(\olapla^{s+l-2}\tau(\p), \onabla_{e_i}\olapla^{s-l-1}\tau(\p) ) \onabla_{e_i}W.\\
\end{align*}
\end{thm}

\begin{proof}
By $(\ref{s})$, we have
\begin{equation}\label{s1}
\begin{split}
\frac{1}{2}\frac{\partial ^2}{\partial r \partial t}& E_{2s}(F)\\
=&\int_M \langle \onabla_{\frac{\partial}{\partial r}} dF(\fpt), -\olapla^{2s-1}\tau(F)
+R^N(\olapla^{2s-2}\tau(F),dF(e_i ) )dF(e_i )\\
&\hspace{10pt}+\sum^{s-1}_{l=1} \{
R^N(\onabla_{e_i}\olapla^{s+l-2}\tau(F) ,\olapla^{s-l-1}\tau(F)  )dF(e_i)\\
&\hspace{10pt}-R^N(\olapla^{s+l-2}\tau(F) ,\onabla_{e_i}\olapla^{s-l-1}\tau(F) )dF(e_i )\} \rangle v_g.\\
&+\int_M \langle 
 F(\fpt), \onabla_{\frac{\partial}{\partial r}} \{-\olapla^{2s-1}\tau(F)
+R^N(\olapla^{2s-2}\tau(F),dF(e_i ) )dF(e_i )\\
&\hspace{10pt}+\sum^{s-1}_{l=1} \{
R^N(\onabla_{e_i}\olapla^{s+l-2}\tau(F) ,\olapla^{s-l-1}\tau(F)  )dF(e_i)\\
&\hspace{10pt}-R^N(\olapla^{s+l-2}\tau(F) ,\onabla_{e_i}\olapla^{s-l-1}\tau(F) )dF(e_i )\}\} \rangle v_g.
\end{split}
\end{equation}
Then, putting t=0, the first term of $(\ref{s1})$ vanishes. Thus, we calculate the second term of $(\ref{s1})$

Using Lemma $\ref{laps-1}$, we have 
\begin{align*}
\onabla_{\frac{\partial}{\partial r}}\olapla^{2s-1}\tau(F)|_{t=0}=I_{2s}.
\end{align*}

\begin{align*}
\onabla_{\frac{\partial}{\partial r}}&R^N(\olapla^{2s-2}\tau(F),dF(e_i))dF(e_i)\\
=&(\nabla^N_{dF(\frac{\partial}{\partial r})}R^N)(\olapla^{2s-2},dF(e_i))dF(e_i)\\
&+R^N(\onabla_{\frac{\partial}{\partial r}}\olapla^{2s-2},dF(e_i))dF(e_i))\\
&+R^N(\olapla^{2s-2},\onabla_{\frac{\partial}{\partial r}} dF(e_i))dF(e_i))\\
&+R^N(\olapla^{2s-2}, dF(e_i))\onabla_{\frac{\partial}{\partial r}}dF(e_i)),\\
\end{align*}
Using second Bianch's identity, Lemma $\ref{laps-1}$, we have
\begin{align*}
\onabla_{\frac{\partial}{\partial r}}&R^N(\olapla^{2s-2}\tau(F),dF(e_i))dF(e_i)|_{t=0}
=II_{2s}.
\end{align*}

\begin{align*}
\onabla_{\frac{\partial}{\partial r}}&R^N(\onabla_{e_i}\olapla^{s+l-2}\tau(F),\olapla^{s-l-1}\tau(F))dF(e_i)\\
=&(\nabla^N_{dF(\frac{\partial}{\partial r})}R^N)(\onabla_{e_i}\olapla^{s+l-2}\tau(F),\olapla^{s-l-1}\tau(F))dF(e_i)\\
&+R^N(\onabla_{\frac{\partial}{\partial r}}\onabla_{e_i}\olapla^{s+l-2}\tau(F),\olapla^{s-l-1}\tau(F))dF(e_i))\\
&+R^N(\onabla_{e_i}\olapla^{s+l-2}\tau(F),\onabla_{\frac{\partial}{\partial r}} \olapla^{s-l-1}\tau(F))dF(e_i))\\
&+R^N(\onabla_{e_i}\olapla^{s+l-2}\tau(F), \olapla^{s-l-1}\tau(F))\onabla_{\frac{\partial}{\partial r}}dF(e_i)),\\
\end{align*}
Using second Bianch's identity, Lemma $\ref{laps-1}$ and Lemma $\ref{nablaps-1}$, we have
$$\onabla_{\frac{\partial}{\partial r}}R^N(\onabla_{e_i}\olapla^{s+l-2}\tau(F),\olapla^{s-l-1}\tau(F))dF(e_i)
=III_{2s}.$$

\begin{align*}
\onabla_{\frac{\partial}{\partial r}}&R^N(\olapla^{s+l-2}\tau(F),\onabla_{e_i}\olapla^{s-l-1}\tau(F))dF(e_i)\\
=&(\nabla^N_{dF(\frac{\partial}{\partial r})}R^N)(\olapla^{s+l-2}\tau(F),\onabla_{e_i}\olapla^{s-l-1}\tau(F))dF(e_i)\\
&+R^N(\onabla_{\frac{\partial}{\partial r}}\olapla^{s+l-2}\tau(F),\onabla_{e_i}\olapla^{s-l-1}\tau(F))dF(e_i))\\
&+R^N(\olapla^{s+l-2}\tau(F),\onabla_{\frac{\partial}{\partial r}} \onabla_{e_i}\olapla^{s-l-1}\tau(F))dF(e_i))\\
&+R^N(\olapla^{s+l-2}\tau(F), \onabla_{e_i}\olapla^{s-l-1}\tau(F))\onabla_{\frac{\partial}{\partial r}}dF(e_i)),\\
\end{align*}
Using second Bianch's identity, Lemma $\ref{laps-1}$ and Lemma $\ref{nablaps-1}$, we have
$$\onabla_{\frac{\partial}{\partial r}}R^N(\olapla^{s+l-2}\tau(F),\onabla_{e_i}\olapla^{s-l-1}\tau(F))dF(e_i)
=IV_{2s}.$$

\end{proof}

\vspace{10pt}

\begin{thm}
Let $\p: (M,g)\rightarrow (N,h)$ be a $(2s+1)$-harmonic map $(s=0,1,\cdots )$. Then, the 
 {\em Hessian} of the $(2s+1)$-energy $E_{2s+1}$ at $\p$ is given by 
\begin{align*}
H(E_{2s+1})_{\p}(V,W)=\int_M\langle V, J_{2s+1}(W)\rangle v_g,
\end{align*}  
where,
\begin{align*}
J_{2s+1}(W)=-I_{2s+1}+II_{2s+1}+III_{2s+1}+IV_{2s+1}+V_{2s+1},
\end{align*}
where,
\begin{align*}
I_{2s+1}=&-\olapla^{2s+1} W
+\olapla^{2s}R^N(W,d\p(e_j) )d\p(e_j)\\
&+\sum^{2s}_{l=1}\olapla^{l-1}
\{
-\onabla_{e_j}R^N(W,d\p(e_j) )\olapla^{2s-l}\tau(\p) \\
&\hspace{52pt}-R^N(W,d\p(e_j) )\onabla_{e_j}\olapla^{2s-l}\tau(\p)\\
&\hspace{52pt}+R^N(W,d\p(\nabla_{e_j}e_j )\olapla^{2s-l}\tau(\p)
\},
\end{align*}

\begin{align*}
II_{2s+1}=&-(\nabla^N_{\olapla^{2s-1}\tau(\p)}R^N)(d\p(e_i),W)d\p(e_i)\\
&-(\nabla^N_{d\p(e_i)}R^N)(W,\olapla^{2s-1}\tau(\p))d\p(e_i)\\
&\hspace{30pt}+R^N(-\olapla^{2s} W
+\olapla^{2s-1}R^N(W,d\p(e_j) )d\p(e_j)\\
&+\sum^{2s-1}_{l_2=1}\{\olapla^{l_2-1}
\{
-\onabla_{e_j}R^N(W,d\p(e_j) )\olapla^{2s-l_2-1}\tau(\p) \\
&\hspace{52pt}-R^N(W,d\p(e_j) )\onabla_{e_j}\olapla^{2s-l_2-1}\tau(\p)\\
&\hspace{52pt}+R^N(W,d\p(\nabla_{e_j}e_j )\olapla^{2s-l_2-1}\tau(\p)
\}\}, d\p(e_i) )d\p(e_i)\\
&+R^N(\olapla^{2s-1}\tau(\p), \onabla_{e_i}W)d\p(e_i)\\
&+R^N(\olapla^{2s-1}\tau(\p), d\p(e_i) ) \onabla_{e_i}W,
\end{align*}

\begin{align*}
III_{2s+1}=&
-(\nabla^N_{\onabla_{e_i}\olapla^{s+l-1}\tau(\p)}R^N)(\olapla^{s-l-1}\tau(\p),W)d\p(e_i)\\
&-(\nabla^N_{\olapla^{s-l-1}}R^N)(W,\onabla_{e_i}\olapla^{s+l-1}\tau(\p))d\p(e_i)\\
&+R^N(-\onabla_{e_i}\olapla^{s+l} W
+\onabla_{e_i}\olapla^{s+l-1}R^N(W,d\p(e_j) )d\p(e_j)\\
&\hspace{30pt}+\sum^{s+l-1}_{l_2=1}\{\onabla_{e_i}\olapla^{l_2-1}
\{
-\onabla_{e_j}R^N(W,d\p(e_j) )\olapla^{s+l-1-l_2}\tau(\p) \\
&\hspace{52pt}-R^N(W,d\p(e_j) )\onabla_{e_j}\olapla^{s+l-1-l_2}\tau(\p)\\
&\hspace{52pt}+R^N(W,d\p(\nabla_{e_j}e_j )\olapla^{s+l-1-l_2}\tau(\p)
\}\\
&+R^N(W,d\p(e_i))\olapla^{s+l-1}\tau(\p)\},\olapla^{s-l-1}\tau(\p))d\p(e_i)\\
&+R^N(\onabla_{e_i}\olapla^{s+l-1}\tau(\p), 
\{-\olapla^{s-l} W
+\olapla^{s-l-1}R^N(W,d\p(e_j) )d\p(e_j)\\
&\hspace{30pt}+\sum^{s-l-1}_{l_2=1}\{\olapla^{l_2-1}
\{
-\onabla_{e_j}R^N(W,d\p(e_j) )\olapla^{s-l-1-l_2}\tau(\p) \\
&\hspace{52pt}-R^N(W,d\p(e_j) )\onabla_{e_j}\olapla^{s-l-1-l_2}\tau(\p)\\
&\hspace{52pt}+R^N(W,d\p(\nabla_{e_j}e_j )\olapla^{s-l-1-l_2}\tau(\p)
\}
\})d\p(e_i)\\
&+R^N(\onabla_{e_i}\olapla^{s+l-1}\tau(\p), \olapla^{s-l-1}\tau(\p) ) \onabla_{e_i}W,\\
\end{align*}

\begin{align*}
IV_{2s}=&
-(\nabla^N_{\olapla^{s+l-1}\tau(\p)}R^N)(\onabla_{e_i}\olapla^{s-l-1}\tau(\p),W)d\p(e_i)\\
&-(\nabla^N_{\onabla_{e_i}\olapla^{s-l-1}}R^N)(W,\olapla^{s+l-1}\tau(\p))d\p(e_i)\\
&+R^N(-\olapla^{s+l} W
+\olapla^{s+l-1}R^N(W,d\p(e_j) )d\p(e_j)\\
&\hspace{30pt}+\sum^{s+l-1}_{l_2=1}\{\olapla^{l_2-1}
\{
-\onabla_{e_j}R^N(W,d\p(e_j) )\olapla^{s+l-1-l_2}\tau(\p) \\
&\hspace{52pt}-R^N(W,d\p(e_j) )\onabla_{e_j}\olapla^{s+l-1-l_2}\tau(\p)\\
&\hspace{52pt}+R^N(W,d\p(\nabla_{e_j}e_j )\olapla^{s+l-1-l_2}\tau(\p)
\},\onabla_{e_i}\olapla^{s-l-1}\tau(\p) )d\p(e_i)\\
&+R^N(\olapla^{s+l-1}\tau(\p), 
\{-\onabla_{e_i}\olapla^{s-l} W
+\onabla_{e_i}\olapla^{s-l-1}R^N(W,d\p(e_j) )d\p(e_j)\\
&\hspace{30pt}+\sum^{s-l-1}_{l_2=1}\{\onabla_{e_i}\olapla^{l_2-1}
\{
-\onabla_{e_j}R^N(W,d\p(e_j) )\olapla^{s-l-1-l_2}\tau(\p) \\
&\hspace{52pt}-R^N(W,d\p(e_j) )\onabla_{e_j}\olapla^{s-l-1-l_2}\tau(\p)\\
&\hspace{52pt}+R^N(W,d\p(\nabla_{e_j}e_j )\olapla^{s-l-1-l_2}\tau(\p)
\}\\
&+R^N(W,d\p(e_i))\olapla^{s-l-1}\tau(\p)\})d\p(e_i)\\
&+R^N(\olapla^{s+l-1}\tau(\p), \onabla_{e_i}\olapla^{s-l-1}\tau(\p) ) \onabla_{e_i}W,
\end{align*}

\begin{align*}
V_{2s+1}=&
-(\nabla^N_{\onabla_{e_i}\olapla^{s-l}\tau(\p)}R^N)(\olapla^{s-1}\tau(\p),W)d\p(e_i)\\
&-(\nabla^N_{\olapla^{s-1}\tau(\p)}R^N)(W,\onabla_{e_i}\olapla^{s-1}\tau(\p))d\p(e_i)\\
&+R^N(-\onabla_{e_i}\olapla^{s} W
+\onabla_{e_i}\olapla^{s-1}R^N(W,d\p(e_j) )d\p(e_j)\\
&\hspace{30pt}+\sum^{s-1}_{l_2=1}\{\onabla_{e_i}\olapla^{l_2-1}
\{
-\onabla_{e_j}R^N(W,d\p(e_j) )\olapla^{s-l_2-1}\tau(\p) \\
&\hspace{115pt}-R^N(W,d\p(e_j) )\onabla_{e_j}\olapla^{s-l_2-1}\tau(\p)\\
&\hspace{115pt}+R^N(W,d\p(\nabla_{e_j}e_j )\olapla^{s-l_2-1}\tau(\p)
\}\\
&+R^N(W,d\p(e_i))\olapla^{s-1}\tau(\p)\},\olapla^{s-1}\tau(\p))d\p(e_i)\\
&+R^N(\onabla_{e_i}\olapla^{s-1}\tau(\p), 
\{-\olapla^{s} W
+\olapla^{s-1}R^N(W,d\p(e_j) )d\p(e_j)\\
&\hspace{30pt}+\sum^{s-1}_{l_2=1}\{\olapla^{l_2-1}
\{
-\onabla_{e_j}R^N(W,d\p(e_j) )\olapla^{s-l_2-1}\tau(\p) \\
&\hspace{52pt}-R^N(W,d\p(e_j) )\onabla_{e_j}\olapla^{s-l_2-1}\tau(\p)\\
&\hspace{52pt}+R^N(W,d\p(\nabla_{e_j}e_j )\olapla^{s-l_2-1}\tau(\p)
\}
\})d\p(e_i)\\
&+R^N(\onabla_{e_i}\olapla^{s-1}\tau(\p), \olapla^{s-1}\tau(\p) ) \onabla_{e_i}W.
\end{align*}
\end{thm}

\begin{proof}
By $(\ref{s+1})$, we have
\begin{equation}\label{s2}
\begin{split}
\frac{1}{2}\frac{\partial ^2}{\partial r \partial t} &E_{2s+1}(F)\\
=&\int_M \langle \onabla_{\frac{\partial}{\partial r}} dF(\fpt), -\olapla^{2s}\tau(F)
+R^N(\olapla^{2s-1}\tau(F),dF(e_j ) )dF(e_j )\\
&\hspace{10pt}+\sum^{s-1}_{l=1} \{
R^N(\onabla_{e_j}\olapla^{s+l-1}\tau(F) ,\olapla^{s-l-1}\tau(F)  )dF(e_j)\\
&\hspace{30pt}-R^N(\olapla^{s+l-1}\tau(F) ,\onabla_{e_j}\olapla^{s-l-1}\tau(F) )dF(e_j )\}\\
&\hspace{10pt}+R^N(\onabla_{e_i}\olapla^{s-1}\tau(F),\olapla^{s-1}\tau(F))dF(e_i)
 \rangle v_g.\\
&+\int_M \langle 
 F(\fpt), \onabla_{\frac{\partial}{\partial r}} \{ -\olapla^{2s}\tau(F)
+R^N(\olapla^{2s-1}\tau(F),dF(e_j ) )dF(e_j )\\
&\hspace{10pt}+\sum^{s-1}_{l=1} \{
R^N(\onabla_{e_j}\olapla^{s+l-1}\tau(F) ,\olapla^{s-l-1}\tau(F)  )dF(e_j)\\
&\hspace{30pt}-R^N(\olapla^{s+l-1}\tau(F) ,\onabla_{e_j}\olapla^{s-l-1}\tau(F) )dF(e_j )\}\\
&\hspace{10pt}+R^N(\onabla_{e_i}\olapla^{s-1}\tau(F),\olapla^{s-1}\tau(F))dF(e_i)\}\} \rangle v_g.
\end{split}
\end{equation}
Then, putting t=0, the first term of $(\ref{s2})$ vanishes. Thus, we calculate the second term of $(\ref{s2})$

Using Lemma $\ref{laps-1}$, we have 
\begin{align*}
\onabla_{\frac{\partial}{\partial r}}\olapla^{2s}\tau(F)|_{t=0}=I_{2s+1}.
\end{align*}

\begin{align*}
\onabla_{\frac{\partial}{\partial r}}&R^N(\olapla^{2s-1}\tau(F),dF(e_j))dF(e_j)\\
=&(\nabla^N_{dF(\frac{\partial}{\partial r})}R^N)(\olapla^{2s-1},dF(e_j))dF(e_j)\\
&+R^N(\onabla_{\frac{\partial}{\partial r}}\olapla^{2s-1},dF(e_j))dF(e_j))\\
&+R^N(\olapla^{2s-1},\onabla_{\frac{\partial}{\partial r}} dF(e_j))dF(e_j))\\
&+R^N(\olapla^{2s-1}, dF(e_j))\onabla_{\frac{\partial}{\partial r}}dF(e_j)),\\
\end{align*}
Using second Bianch's identity, Lemma $\ref{laps-1}$, we have
\begin{align*}
\onabla_{\frac{\partial}{\partial r}}&R^N(\olapla^{2s-1}\tau(F),dF(e_j))dF(e_j)|_{t=0}
=II_{2s+1}.
\end{align*}

\begin{align*}
\onabla_{\frac{\partial}{\partial r}}&R^N(\onabla_{e_i}\olapla^{s+l-1}\tau(F),\olapla^{s-l-1}\tau(F))dF(e_j)\\
=&(\nabla^N_{dF(\frac{\partial}{\partial r})}R^N)(\onabla_{e_i}\olapla^{s+l-1}\tau(F),\olapla^{s-l-1}\tau(F))dF(e_j)\\
&+R^N(\onabla_{\frac{\partial}{\partial r}}\onabla_{e_i}\olapla^{s+l-1}\tau(F),\olapla^{s-l-1}\tau(F))dF(e_j))\\
&+R^N(\onabla_{e_i}\olapla^{s+l-1}\tau(F),\onabla_{\frac{\partial}{\partial r}} \olapla^{s-l-1}\tau(F))dF(e_j))\\
&+R^N(\onabla_{e_i}\olapla^{s+l-1}\tau(F), \olapla^{s-l-1}\tau(F))\onabla_{\frac{\partial}{\partial r}}dF(e_j)),\\
\end{align*}
Using second Bianch's identity, Lemma $\ref{laps-1}$ and Lemma $\ref{nablaps-1}$, we have
$$\onabla_{\frac{\partial}{\partial r}}R^N(\onabla_{e_i}\olapla^{s+l-1}\tau(F),\olapla^{s-l-1}\tau(F))dF(e_j)
=III_{2s+1}.$$

\begin{align*}
\onabla_{\frac{\partial}{\partial r}}&R^N(\olapla^{s+l-1}\tau(F),\onabla_{e_i}\olapla^{s-l-1}\tau(F))dF(e_j)\\
=&(\nabla^N_{dF(\frac{\partial}{\partial r})}R^N)(\olapla^{s+l-1}\tau(F),\onabla_{e_i}\olapla^{s-l-1}\tau(F))dF(e_j)\\
&+R^N(\onabla_{\frac{\partial}{\partial r}}\olapla^{s+l-1}\tau(F),\onabla_{e_i}\olapla^{s-l-1}\tau(F))dF(e_j))\\
&+R^N(\olapla^{s+l-1}\tau(F),\onabla_{\frac{\partial}{\partial r}} \onabla_{e_i}\olapla^{s-l-1}\tau(F))dF(e_j))\\
&+R^N(\olapla^{s+l-1}\tau(F), \onabla_{e_i}\olapla^{s-l-1}\tau(F))\onabla_{\frac{\partial}{\partial r}}dF(e_j)),\\
\end{align*}
Using second Bianch's identity, Lemma $\ref{laps-1}$ and Lemma $\ref{nablaps-1}$, we have
$$\onabla_{\frac{\partial}{\partial r}}R^N(\olapla^{s+l-1}\tau(F),\onabla_{e_i}\olapla^{s-l-1}\tau(F))dF(e_j)
=IV_{2s+1}.$$

\begin{align*}
\onabla_{\frac{\partial}{\partial r}}&R^N(\onabla_{e_i}\olapla^{s-1}\tau(F),\olapla^{s-1}\tau(F))dF(e_j)\\
=&(\nabla^N_{dF(\frac{\partial}{\partial r})}R^N)(\onabla_{e_i}\olapla^{s-1}\tau(F),\olapla^{s-1}\tau(F))dF(e_j)\\
&+R^N(\onabla_{\frac{\partial}{\partial r}}\onabla_{e_i}\olapla^{s-1}\tau(F),\olapla^{s-1}\tau(F))dF(e_j))\\
&+R^N(\onabla_{e_i}\olapla^{s-1}\tau(F),\onabla_{\frac{\partial}{\partial r}} \olapla^{s-1}\tau(F))dF(e_j))\\
&+R^N(\onabla_{e_i}\olapla^{s-1}\tau(F), \olapla^{s-1}\tau(F))\onabla_{\frac{\partial}{\partial r}}dF(e_j)),\\
\end{align*}
Using second Bianch's identity, Lemma $\ref{laps-1}$ and Lemma $\ref{nablaps-1}$, we have
$$\onabla_{\frac{\partial}{\partial r}}R^N(\onabla_{e_i}\olapla^{s-1}\tau(F),\olapla^{s-1}\tau(F))dF(e_j)
=V_{2s+1}.$$

\end{proof}

\begin{defn}
Assume that $\p:(M,g)\rightarrow (N,h)$ is a $k$-harmonic map. 
 The operator $J_{k}$ on $\Gamma (\p^{-1}TN)$ is the 2$k$th order self-adjoint elliptic differential operator, 
 so that it has a spectrum consisting of discrete eigenvalues 
$\lambda_{1} <\lambda_2 <\cdots <\lambda_t<\cdots$ with their finite multiplicities. 
 Denote by $E^k_{\lambda_1},E^k_{\lambda_2}, \cdots, $ the corresponding eigenspaces in $\Gamma(\p^{-1}TN)$. 
 Then, the definitions of {\em k-index} and {\em k-nullity} are given by
$${\rm index}_{k}(\p)={\rm dim}(\oplus_{\lambda <0}E^k_{\lambda}),$$ 
$${\rm nullity}_{k}(\p)={\rm dim}E^k_0.$$
And $\p$ is {\em weakly stable } if ${\rm index}=0 $,
i.e., $H(E_{k})_{\p}(V,V)\geq 0,$
 for all $V \in \Gamma (\p^{-1}TN)$.
 $\p$ is {\em unstable} if it is not {\em weakly stable}.
\end{defn}

\vspace{10pt}

\begin{prop}
Any harmonic map is weakly stable $k$-harmonic map.
\end{prop}

\begin{proof}
Case1. $k=2s,  (s=1,2, \cdots) .$

By assumption we have
\begin{equation}\label{hess}
\begin{split}
H_{2s}(V,V)
=&\int_M \langle V, -(-\olapla^{2s}V+\olapla^{2s-1}R^N(V,d\p(e_j))d\p(e_j))\\
&\hspace{29pt}+R^N(-\olapla^{2s-1}V+\olapla^{2s-2}R^N(V,d\p(e_j))d\p(e_j), d\p(e_i))d\p(e_i)\rangle v_g\\
=&\int_M \langle -V, -\olapla^{2s}V+\olapla^{2s-1}R^N(V,d\p(e_j))d\p(e_j) \rangle v_g\\
&+\int_M \langle R^N(V,d\p(e_i))d\p(e_i),-\olapla^{2s-1}V+\olapla^{2s-2}R^N(V,d\p(e_j))d\p(e_j) \rangle v_g\\
=&\int_M||-\olapla^{s}V+\olapla^{s-1}R^N(V,d\p(e_j))d\p(e_i)||^2 v_g \geq 0.
\end{split}
\end{equation}

Case 2. $k=2s+1,  (s=0,1,2, \cdots) .$

By assumption we have 
\begin{equation}\label{hess+1}
\begin{split}
H_{\p}(V,V)
=&\int_M\langle V, -(-\olapla^{2s+1}V+\olapla^{2s}R^N(V,d\p(e_j))d\p(e_j) ) )\\
&\hspace{29pt}+R^N(-\olapla^{2s}V+\olapla^{2s-1}R^N(V,d\p(e_j))d\p(e_j), d\p(e_i))d\p(e_i)\rangle v_g\\
=&\int_M \langle -V, -\olapla^{2s+1}V+\olapla^{2s}R^N(V,d\p(e_j))d\p(e_j) \rangle v_g\\
&+\int_M \langle R^N(V,d\p(e_i))d\p(e_i),-\olapla^{2s}V+\olapla^{2s-1}R^N(V,d\p(e_j))d\p(e_j) \rangle v_g\\
=&\int_M||\onabla(-\olapla^{s}V+\olapla^{s-1}R^N(V,d\p(e_j))d\p(e_j))||^2 v_g \geq 0.
\end{split}
\end{equation} 

\end{proof}

\begin{cor}
Assume that $\p: (M,g)\rightarrow (N,h)$ is a harmonic map. 
Then, 
$${\rm nullity}_{k}(\p)=\{V\in \Gamma (\p^{-1}TN); J(\olapla^{k-2}J(V))=0\}.$$ 
\end{cor}

\begin{proof}
If $\p$ is harmonic map, then, $\tau(\p)=0.$
Thus we have 
$$J_{k}(V)=J(\olapla^{k-2}J(V)),$$
for all $V\in \Gamma(\p^{-1}TN).$
Therefore, we have the corollary.
\end{proof}


\section{The $k$-harmonic maps into the product spaces}\label{product}

In this section, we describe the necessary and sufficient condition of 
$k$-harmonic maps into the product spaces. 
 First, let us recall the result of Y.-L. Ou \cite{ylo1}.
\begin{thm}[\cite{ylo1}]\label{thm1 of ylo1}
Let $\varphi :(M,g)\rightarrow (N_1,h_1)$ and 
$\psi :(M,g)\rightarrow (N_2,h_2)$ be two maps. Then, the map 
$\phi:(M,g)\rightarrow (N_1\times N_2,h_1\times h_2)$ with 
$\phi(x)=(\varphi(x),\psi(x))$ is $2$-harmonic if and only if the both map 
$\varphi$ or $\psi$ are $2$-harmonic. Furthermore, if one of $\varphi$ or $\psi$ is $2$-harmonic
and the other is a proper $2$-harmonic map, then $\phi$ is a proper $2$-harmonic map.
\end{thm}
We generalize Theorem \ref{thm1 of ylo1} for $k$-harmonic maps. Namely, we have
 the following theorem which is useful to construct examples the $k$-harmonic maps.
\begin{thm}\label{ps}
Let $\varphi :(M,g)\rightarrow (N_1,h_1)$ and 
$\psi :(M,g)\rightarrow (N_2,h_2)$ be two maps. Then, the map 
$\phi:(M,g)\rightarrow (N_1\times N_2,h_1\times h_2)$ with 
$\phi(x)=(\varphi(x),\psi(x))$ is $k$-harmonic if and only if the both map 
$\varphi$ or $\psi$ are $k$-harmonic. Furthermore, if one of $\varphi$ or $\psi$ is harmonic 
and the other is a proper $k$-harmonic map, then $\phi$ is a proper $k$-harmonic map.
\end{thm}

\begin{proof}
It is easily seen that 
\begin{equation}
\begin{split}
d\p(X)=d\varphi (X)+d\phi (X),\ \ \ \forall X\in \Gamma (TM).
\end{split}
\end{equation}

It follows that
\begin{align}
\nabla^{\p}_{X}d\p (Y)= \nabla^{\p}_{X}d\varphi (Y)+\nabla^{\p}_{X}d\psi (Y),\ \ X,Y\in \Gamma(TM).
\end{align}
where,$\nabla^{\p}$ is given by
$\nabla^{\p}_{X}=\nabla^N_{d\p (X)},$ $\forall X\in \Gamma (TM)$.

Let $\{e_i\}_{i=1}^m$ be a local orthonaormal frame on $(M,g)$ and $Y=Y^ie_i$, then 
$d\varphi (Y)=Y^i\varphi^{\alpha}_i(E_{\alpha }\varphi)$, 
for some function $\varphi ^{\alpha}$ defined locally on $M$.
A straight forward computation yields 
$\nabla^{\p}_{X}d\varphi (Y)=\nabla^{\varphi}_{X}d\varphi (Y)$,
 $\tau(\phi)=\tau(\varphi)+\tau(\psi),$
 $\overline \triangle_{\phi}\tau(\phi)=\overline \triangle_{\varphi}\tau(\varphi)+\overline \triangle_{\psi}\tau(\psi).$
And we notice that 
$\overline \triangle_{\varphi}\tau(\varphi)$ is tangent to $N_1$, $\overline \triangle_{\psi}\tau(\psi)$ is tangent to $N_2$.
So we have 
\begin{equation}
\begin{split}
\overline \triangle_{\phi} (\overline \triangle_{\phi}\tau(\phi))
&=\overline \triangle_{\phi}(\overline \triangle_{\varphi}\tau(\varphi)+\overline \triangle_{\psi}\tau(\psi))\\
&=\overline \triangle_{\varphi}(\overline \triangle_{\varphi}\tau(\varphi))
+\overline \triangle_{\psi}(\overline \triangle_{\psi}\tau(\psi)).
\end{split}
\end{equation}
Simillary,
\begin{align*}
\overline \triangle^{t}_{\phi}\tau(\phi)
=\overline \triangle^{t}_{\varphi}\tau(\varphi)
+\overline \triangle^{t}_{\psi}\tau(\psi).
\end{align*}
for all $t=0,1,2,\cdots.$

We use the property of the curvature of the product manifold to have
\begin{align*}
&R^{N_1\times N_2}(\overline \triangle ^{t}_{\phi}\tau(\phi),d\phi(e_i))d\phi(e_i)\\
=&R^{N_1}(\overline \triangle^{t}_{\varphi}\tau(\varphi),d\varphi(e_i))d\varphi(e_i)
+R^{N_2}(\overline \triangle^{t}_{\psi}\tau(\psi),d\psi(e_i))d\psi(e_i).
\end{align*}
Simillary we have

\begin{align*}
&R^{N_1\times N_2}(\overline \triangle ^{s}_{\phi}\tau(\phi),\overline \triangle ^{t}_{\phi}\tau(\phi))\\
=&R^{N_1}(\overline \triangle^{s}_{\varphi}\tau(\varphi),\overline \triangle^{t}_{\varphi}\tau(\varphi))
+R^{N_2}(\overline \triangle^{s}_{\psi}\tau(\psi),\overline \triangle^{t}_{\psi}\tau(\psi)),
\end{align*}

\begin{align*}
&R^{N_1\times N_2}(\nabla^{\p}_{d\p(X)}\overline \triangle ^{s}_{\phi}\tau(\phi),\overline \triangle ^{t}_{\phi}\tau(\phi))\\
=&R^{N_1}(\nabla^{\varphi}_{d\varphi (X)}\overline \triangle^{s}_{\varphi}\tau(\varphi),\overline \triangle^{t}_{\varphi}\tau(\varphi))
+R^{N_2}(\nabla^{\psi}_{d\psi (X)}\overline \triangle^{s}_{\psi}\tau(\psi),\overline \triangle^{t}_{\psi}\tau(\psi)),
\end{align*}

\begin{align*}
&R^{N_1\times N_2}(\overline \triangle ^{s}_{\phi}\tau(\phi),\nabla^{\p}_{d\p(X)}\overline \triangle ^{t}_{\phi}\tau(\phi))\\
=&R^{N_1}(\overline \triangle^{s}_{\varphi}\tau(\varphi),\nabla^{\varphi}_{d\varphi (X)}\overline \triangle^{t}_{\varphi}\tau(\varphi))
+R^{N_2}(\overline \triangle^{s}_{\psi}\tau(\psi),\nabla^{\psi}_{d\psi (X)}\overline \triangle^{t}_{\psi}\tau(\psi)),
\end{align*}
for all $t,s=0,1,2,\cdots $, and for all $X\in \Gamma(TM)$.

Thus using Theorem $\ref{2s-harmonic}$, $\ref{2s+1-harmonic}$, we have the theorem.
\end{proof}

\vspace{10pt}
The following corollary generalizes Corollary 3.4 in \cite{ylo1}.

\begin{cor}
Let $\psi:(M,g)\rightarrow (N,h)$ be a smooth map. Then, the graph 
$\p:(M,g)\rightarrow (M\times N,g\times h)$ with $\p(x)=(x,\psi(x))$ is a $k$-harmonic map 
 if and only if the map $\psi:(M,g)\rightarrow (N,h)$ is a $k$-harmonic map.
 Furthermore, if $\psi$ is proper $k$-harmonic, then so is the graph.  
\end{cor}

\begin{proof}
This follows from Theorem \ref{ps} with $\varphi:(M,g)\rightarrow (N,h)$ being 
 identity map which is harmonic.
\end{proof}




\section{$k$-harmonic curves into a Riemannian manifold with constant sectional curvature}\label{constant}

As well known harmonic map always biharmonic map. And by Corollary $\ref{harmonick-harmonic}$, harmonic map is always $k$-harmonic map. In this section we consider the next problem.
\begin{prob}\label{probst}
biharmonic map is always $k$-harmonic map $(k=2,3,\cdots )$ ?
More generally, for $s<t$, 
$s$-harmonic map is always $t$-harmonic map ?
\end{prob}





First, the Frenet-frame is given as follows:
\begin{equation}
\begin{cases}
&\gamma '=T\ ,\ \nabla^N _{\gamma '}T=\kappa N\ ,\ \nabla^N_{\gamma '}N=-\kappa T,\\
&\langle T,N \rangle =0\ ,\ 
\langle T,T \rangle =1\ ,\ 
\langle N,N \rangle =1,
\end{cases}
\end{equation}
where $\kappa$ is the geodesic curvature and $\langle\cdot,\cdot\rangle=h$ the Riemannian metric on $N$.  
 Then, we have the following.

\vspace{10pt}

\begin{prop}\label{8}
Let $\gamma:I\rightarrow (N^2,\langle \cdot, \cdot \rangle)$ be a smooth curve parametrized by arc length 
 from an open interval of $\mathbb{R}$ into a Riemannian manifold $(N^2,\langle \cdot, \cdot \rangle)$ with constant sectional curvature $K$.
Then, $\gamma $ is a $3$-harmonic curve if and only if
\begin{equation*}
\begin{cases}
\kappa^{(4)}-12(\kappa')^2-10\kappa^2\kappa''+\kappa^5-3\kappa(\kappa')^2
+K(\kappa''-2\kappa^3)=0,\\
\kappa\kappa^{(3)}-2\kappa^3\kappa'+2\kappa'\kappa''=0,
\end{cases}
\end{equation*}
where $\kappa $ is the geodesic curvature of $\gamma$.
\end{prop}


\begin{proof}
We calculate
 $(\nabla^N_{\gamma'}\nabla^N_{\gamma'})^2\tau(\gamma)$ as follows.
\begin{equation}\label{nabla22}
\begin{split}
(\nabla^N_{\gamma'}\nabla^N_{\gamma'})^2\tau(\gamma)
&=(\kappa^{(4)}-12(\kappa')^2-10\kappa^2\kappa''+\kappa^5-3\kappa(\kappa')^2)N\\
&\ \ +(-5\kappa\kappa^{(3)}+10\kappa^3\kappa'-10\kappa'\kappa'')T.
\end{split}
\end{equation}
Therefore, $\gamma$ is $3$-harmonic if and only if
\begin{equation}
\begin{split}
&(\kappa^{(4)}-12(\kappa')^2-10\kappa^2\kappa''+\kappa^5-3\kappa(\kappa')^2+K(\kappa''-2\kappa^3))N\\
&+(-5\kappa\kappa^{(3)}+10\kappa^3\kappa'-10\kappa'\kappa'')T=0.
\end{split}
\end{equation}
So we have Proposition \ref{8}.
\end{proof}

\vspace{10pt}

\begin{cor}
Let $\gamma:I\rightarrow (N^2,\langle \cdot, \cdot \rangle)$ be a $3$-harmonic curve parametrized by arc length 
 from an open interval of $\mathbb{R}$ into a Riemannian manifold $(N^2,\langle \cdot, \cdot \rangle)$ with constant sectional curvature $K\geq 0$.
If geodesic curvature $\kappa$ is constant, $\kappa =\sqrt{2K}$.
\end{cor}

\begin{proof}
We can show this corollary by a direct computation. The proof is omitted.
\end{proof}

\begin{prop}\label{2s-harmonic c}
Let $\gamma:I\rightarrow (N^n,\langle \cdot, \cdot \rangle)$ be a smooth curve parametrized by arc length 
 from an open interval of $\mathbb{R}$ into a Riemannian manifold $(N^n,\langle \cdot, \cdot \rangle)$ with constant sectional curvature $K$.
Then, $\g$ is $2s$-harmonic curve if and only if

\begin{equation}\label{2s-harmonic curve}
\begin{split}
\tau_{2s}(\gamma )
=&(\nabla^N_{\gamma '}\nabla^N_{\gamma '})^{2s-1}\tau(\g)\\
&+K\{(\nabla^N_{\gamma'}\nabla^N_{\gamma'})^{2s-2}\tau(\g)
-\langle \g ', (\nabla^N_{\gamma'}\nabla^N_{\gamma'})^{2s-2}\tau(\g)\rangle \g '\}\\
&-\sum^{s-1}_{l=1}K\{(\langle (\nabla^N_{\gamma'}\nabla^N_{\gamma'})^{s-l-1}\tau(\g),\g' \rangle 
\nabla^N_{\g '}(\nabla^N_{\gamma'}\nabla^N_{\gamma'})^{s+l-2}\tau(\g))\}\\
&\hspace{20pt}-\langle \g ' , \nabla^N_{\g '}(\nabla^N_{\gamma'}\nabla^N_{\gamma'})^{s+l-2}\tau(\g) \rangle 
(\nabla^N_{\gamma'}\nabla^N_{\gamma'})^{s-l-1}\tau(\g))\\
&\hspace{20pt}-(\langle \nabla^N_{\g '}(\nabla^N_{\gamma'}\nabla^N_{\gamma'})^{s-l-1}\tau(\g)) , \g ' \rangle 
(\nabla^N_{\gamma'}\nabla^N_{\gamma'})^{s+l-2}\tau(\g)\\
&\hspace{20pt}+\langle \g ' , (\nabla^N_{\gamma'}\nabla^N_{\gamma'})^{s+l-2}\tau(\g) \rangle 
\nabla^N_{\g'}(\nabla^N_{\gamma'}\nabla^N_{\gamma'})^{s-l-1}\tau(\g)\}=0.
\end{split}
\end{equation}
\end{prop}

\begin{proof}
We only notice that
$$\olapla=-\nabla^N_{\g '}\nabla^N_{\gamma '},$$
$$R^N(V,W)Z=K(\langle W, Z\rangle V -\langle Z, V \rangle W),$$
$$\langle \g', \g' \rangle =1.$$
We get the proposition.
\end{proof}

\vspace{10pt}

Similarly we have

\begin{prop}\label{2s+1-harmonic c}
Let $\gamma:I\rightarrow (N^n,\langle \cdot, \cdot \rangle)$ be a smooth curve parametrized by arc length 
 from an open interval of $\mathbb{R}$ into a Riemannian manifold $(N^n,\langle \cdot, \cdot \rangle)$ with constant sectional curvature $K$.
Then, $\g$ is $(2s+1)$-harmonic curve if and only if
\begin{equation}\label{2s+1-harmonic curve}
\begin{split}
\tau_{2s+1}(\gamma )
=&-(\nabla^N_{\gamma '}\nabla^N_{\gamma '})^{2s}\tau(\g)\\
&-K\{(\nabla^N_{\gamma'}\nabla^N_{\gamma'})^{2s-1}\tau(\g)
-\langle \g ', (\nabla^N_{\gamma'}\nabla^N_{\gamma'})^{2s-1}\tau(\g)\rangle \g '\}\\
&+\sum^{s-1}_{l=1}K\{(\langle (\nabla^N_{\gamma'}\nabla^N_{\gamma'})^{s-l-1}\tau(\g),\g '\rangle 
\nabla^N_{\g '}(\nabla^N_{\gamma'}\nabla^N_{\gamma'})^{s+l-1}\tau(\g))\}\\
&\hspace{20pt}-\langle \g ' , \nabla^N_{\g '}(\nabla^N_{\gamma'}\nabla^N_{\gamma'})^{s+l-1}\tau(\g) \rangle 
(\nabla^N_{\gamma'}\nabla^N_{\gamma'})^{s-l-1}\tau(\g))\\
&\hspace{20pt}-(\langle \nabla^N_{\g '}(\nabla^N_{\gamma'}\nabla^N_{\gamma'})^{s-l-1}\tau(\g)) , \g ' \rangle 
(\nabla^N_{\gamma'}\nabla^N_{\gamma'})^{s+l-1}\tau(\g)\\
&\hspace{20pt}+\langle \g ' , (\nabla^N_{\gamma'}\nabla^N_{\gamma'})^{s+l-1}\tau(\g) \rangle 
\nabla^N_{\g'}(\nabla^N_{\gamma'}\nabla^N_{\gamma'})^{s-l-1}\tau(\g)\}\\
&\hspace{20pt}+(\langle (\nabla^N_{\gamma'}\nabla^N_{\gamma'})^{s-1}\tau(\g),\g '\rangle 
\nabla^N_{\g '}(\nabla^N_{\gamma'}\nabla^N_{\gamma'})^{s-1}\tau(\g))\}\\
&\hspace{20pt}-\langle \g ' , \nabla^N_{\g '}(\nabla^N_{\gamma'}\nabla^N_{\gamma'})^{s-1}\tau(\g) \rangle 
(\nabla^N_{\gamma'}\nabla^N_{\gamma'})^{s-1}\tau(\g))\}=0.
\end{split}
\end{equation}
\end{prop}

\vspace{10pt}

Using these propositions, we show the following propositions. 

\begin{prop}
Let $\gamma:I\rightarrow (N^2,\langle \cdot, \cdot \rangle)$ be a $2s$-harmonic curve $(s=1,2,\cdots )$ parametrized by arc length 
 from an open interval of $\mathbb{R}$ into a Riemannian manifold $(N^2,\langle \cdot, \cdot \rangle)$ with constant sectional curvature $K\geq 0$.
If geodesic curvature $\kappa$ is constant, $\kappa =\sqrt{(2s-1)K}$.
\end{prop}

\begin{proof}
By assumption, for all $t=0,1,2,\cdots $
\begin{align*}
&(\nabla_{\g '}\nabla_{\g'})^t\tau(\g)=(-1)^t\ka^{2t+1} N,\\
&\nabla_{\g'}(\nabla_{\g '}\nabla_{\g'})^t\tau(\g)=-(-1)^t\ka^{2t+2} N.
\end{align*}
Using  these and Proposition $\ref{2s-harmonic c}$, we have
\begin{align*}
\tau_{2s}(\g )
=&-\ka^{4s-1}N+K(\ka^{4s-3}N+2K(s-1)\ka^{4s-3}N)\\
=&\ka^{4s-3}(-\ka^2+(2s-1)K)N=0.
\end{align*}
Therefore we have the proposition.
\end{proof}

\begin{prop}
Let $\gamma:I\rightarrow (N^2,\langle \cdot, \cdot \rangle)$ be a $(2s+1)$-harmonic curve $(s=0,1,2,\cdots )$ parametrized by arc length 
 from an open interval of $\mathbb{R}$ into a Riemannian manifold $(N^2,\langle \cdot, \cdot \rangle)$ with constant sectional curvature $K\geq 0$.
If geodesic curvature $\kappa$ is constant, $\kappa =\sqrt{2sK}$.
\end{prop}

\begin{proof}
By assumption, for all $t=0,1,2\cdots $,
\begin{align*}
&(\nabla_{\g '}\nabla_{\g'})^t\tau(\g)=(-1)^t\ka^{2t+1} N,\\
&\nabla_{\g'}(\nabla_{\g '}\nabla_{\g'})^t\tau(\g)=-(-1)^t\ka^{2t+2} N.
\end{align*}
Using  these and Proposition $\ref{2s+1-harmonic c}$, we have
\begin{align*}
\tau_{2s+1}(\g )
=&-\ka^{4s+1}N+K\ka^{4s-1}N+2K(s-1)\ka^{4s-1}N+K\ka^{4s-1}N\\
=&\ka^{4s-1}\{-\ka^2+2sK\}N=0.
\end{align*}
Therefore we have the proposition.
\end{proof}

\vspace{10pt} 

Therefore, we get the answer of Problem $\ref{probst}$.
For $s<t,$ $s$-harmonic map is not always $t$-harmonic map.

\vspace{10pt} 

Fianlly, we determine that the ODEs of the $3$-harmonic curve equations into a sphere.

\begin{prop}\label{5}
Let $\gamma :I\rightarrow S^n\subset \mathbb{R}^{n+1}$ be a smooth curve parametrized by arc length. 
Then $\gamma$ is $3$-harmonic curve if and only if 
\begin{equation}\label{3harmonic}
\begin{split}
-\gamma ^{(6)}-2\gamma ^{(4)}
-(2g_{13}+3)\gamma ''
+4g_{23}\gamma '
+(1+9g_{24}+8g_{33})\gamma =0,
\end{split}
\end{equation}
where $g_{ij}=g_0(\gamma ^{(i)},\gamma ^{(j)}),(i,j=0,1\dots),$
 and $g_0$ is the standard metric on the Euclidean space $\mathbb{R}^{n+1}$.
\end{prop}


\begin{proof}
$$\nabla ^0_{\gamma '} \gamma '=\sigma (\gamma ',\gamma ')+\nabla _{\gamma '} \gamma ',$$
which yields that
$$\nabla _{\gamma '} \gamma '=\nabla ^0_{\gamma '} \gamma '
+g(\gamma ',\gamma ')\gamma.$$
 Therefore, we have
$\nabla _{\gamma '}\gamma '
=\gamma ''+\gamma .$
 Similarly,
\begin{align}
(\nabla _{\gamma '}\nabla _{\gamma '})(\nabla _{\gamma '}\gamma ')\notag
&=\gamma ^{(4)}+\gamma ''+(g_{13}+1)\gamma .
\end{align}
\begin{equation}\label{nabla2}
\begin{split}
(\nabla _{\gamma '}\nabla _{\gamma '})^2&(\nabla _{\gamma '}\gamma ')\\
&=\gamma ^{(6)}+\gamma ^{(4)}+(g_{13}+1)\gamma ''
+(2g_{23}+3g_{14})\gamma '\\
&\hspace{30pt} +(1+g_{33}+3g_{24}+3g_{15}+g_{22}+3g_{13})\gamma,
\end{split}
\end{equation}

\begin{equation}
\begin{split}
R^N(\olapla\tau(\gamma),\g')\g'=-\g^{(4)}-\g''-(g_{13}+1)\g+g_{14}\g',
\end{split}
\end{equation}

\begin{equation}
\begin{split}
R^N(\nabla_{\g'}\tau(\gamma),\g')\g'=-(g_{13}+1)\g''-(g_{13}+1)\g.
\end{split}
\end{equation}

where, we used
$g_{13}=-g_{22}\ ,\ g_{14}=-3g_{23}\ ,\ g_{15}=-3g_{33}-4g_{24}$.
So we have Proposition {\ref{5}}.
\end{proof}




\begin{thebibliography}{99}
\bibitem{tijihu1}
T. Ichiyama, J. Inoguchi and H. Urakawa,
\textit{Bi-harmonic map and bi-Yang-Mills fields}, Note di Matematica,  {\bf 28} (2009), 233-275.

\bibitem{hu1}
H. Urakawa, \textit{Calculus of variation and harmonic maps}, Transl. Math. Monograph. {\bf 132}, Amer. Math. Soc.


\bibitem{rcsmco1}
R. Caddeo, S. Montaldo and C. Oniciuc, 
\textit{Biharmonic submanifolds in spheres,} 
Israel Journal of Mathmatics, {\bf 130} (2002), 109-123.

\bibitem{ylo1}
Y.-L. Ou, 
\textit{Some constructions of biharmonic maps and Chen's conjecture on biharmonic hypersurfaces},
 arXiv:0912.1141v1 [math.DG] 6 Dec 2009.


\bibitem{jg1}
G. Y. Jiang, 
{\textit 2-harmonic maps and their first and second variational formulas},
Chinese Ann. Math. , 7A (1986), 388-402; the English translation, Note di Matematica, {\bf 28}, (2009), 209-232.

\bibitem{jell1}
J. Eells and L. Lemaire, 
 {\em Selected topics in harmonic maps},
 CBMS, {\bf 50}, Amer. Math. Soc, 1983.
 

\bibitem{pp1}
P. Petersen,
 \textit{Riemannian Geometry},
 Springer Science 2006.
 
\bibitem{ms1}
M. Spivak,
 \textit{A Comprehensive Introduction to Differential Geometry},  {\bf I - IV}, Wilmington: Publish or Perlish, 1979.

\bibitem{ws1}
Wang Shaobo,
\textit{The First Variation Formula For $K$-harmonic mapping},
Journal of jiangxi university, Vol 13, No 1, 1989.



\end{thebibliography}
\end{document}